\documentclass[11pt]{article}

\usepackage{amscd,amsmath, amssymb, fancyhdr, graphicx, url}

\numberwithin{equation}{section}

\newcommand{\version}{version 1.0,\ \ Apr. 2, 2018}

\renewcommand{\leftmark}{{\scriptsize MBM loci in families of hyperk\"ahler manifolds}}

\makeatletter
\def\x@arrow{\DOTSB\Relbar}
\def\xlongrightarrowfill@{\arrowfill@\relbar\relbar\longrightarrow}
\newcommand{\xlongrightarrow}[2][]{%
        \ext@arrow 0099\xlongrightarrowfill@{#1}{#2}}
\makeatother

\def\eqref#1{(\ref{#1})}
\newcommand{\goth}{\mathfrak}

\newcommand{\arrow}{{\:\longrightarrow\:}}
\newcommand{\Z}{{\Bbb Z}}
\def\C{{\Bbb C}}
\def\P{{\Bbb P}}

\newcommand{\R}{{\Bbb R}}
\newcommand{\Q}{{\Bbb Q}}

\def\1{\sqrt{-1}\:}

\newcommand{\cntrct}                
{\hspace{2pt}\raisebox{1pt}{\text{$\lrcorner$}}\hspace{2pt}}

\renewcommand{\bar}{\overline}
\renewcommand{\phi}{\varphi}
\renewcommand{\epsilon}{\varepsilon}
\renewcommand{\geq}{\geqslant}

\renewcommand{\min}{{\operatorname{\sf min}}}
\newcommand{\Teich}{\operatorname{\sf Teich}}

\newcommand{\Comp}{\operatorname{\sf Comp}}

\newcommand{\Per}{\operatorname{\sf Per}}
\newcommand{\Perspace}{\operatorname{{\Bbb P}\sf er}}

\newcommand{\Kah}{\operatorname{Kah}}

\newcommand{\Gr}{\operatorname{Gr}}

\newcommand{\Pos}{\operatorname{Pos}}

\newcommand{\Aut}{\operatorname{Aut}}

\newcommand{\Mon}{\operatorname{\sf Mon}}

\newcommand{\Diff}{\operatorname{\sf Diff}}

\renewcommand{\Re}{\operatorname{Re}}

\newcommand{\proof}{\noindent{\bf Proof:\ }}
\newcommand{\pstep}{{\bf Proof. Step 1:\ }}

\newcounter{Mycounter}[section]
\newcounter{lemma}[section]
\renewcommand{\thelemma}{{Lemma \thesection.\arabic{lemma}}}
\newcommand{\lemma}{%
    \setcounter{lemma}{\value{Mycounter}}
    \refstepcounter{lemma}
    \stepcounter{Mycounter}
    {\noindent \bf \thelemma:\ }}

\newcounter{claim}[section]
\renewcommand{\theclaim}{{Claim \thesection.\arabic{claim}}}
\newcommand{\claim}{%
    \setcounter{claim}{\value{Mycounter}}
    \refstepcounter{claim}
    \stepcounter{Mycounter}
    {\noindent \bf \theclaim:\ }}

\newcounter{sublemma}[section]
\newcounter{corollary}[section]
\renewcommand{\thecorollary}{{Corollary \thesection.\arabic{corollary}}}
\newcommand{\corollary}{%
    \setcounter{corollary}{\value{Mycounter}}
    \refstepcounter{corollary}
    \stepcounter{Mycounter}
    {\noindent \bf \thecorollary:\ }}

\newcounter{theorem}[section]
\renewcommand{\thetheorem}{{Theorem \thesection.\arabic{theorem}}}
\newcommand{\theorem}{%
    \setcounter{theorem}{\value{Mycounter}}
    \refstepcounter{theorem}
    \stepcounter{Mycounter}
    {\noindent \bf \thetheorem:\ }}

\newcounter{conjecture}[section]
\newcounter{proposition}[section]
\renewcommand{\theproposition}
      {{Proposition \thesection.\arabic{proposition}}}
\newcommand{\proposition}{%
    \setcounter{proposition}{\value{Mycounter}}
    \refstepcounter{proposition}
    \stepcounter{Mycounter}
    {\noindent \bf \theproposition:\ }}

\newcounter{definition}[section]
\renewcommand{\thedefinition}
      {{Definition~\thesection.\arabic{definition}}}
\newcommand{\definition}{%
    \setcounter{definition}{\value{Mycounter}}
    \refstepcounter{definition}
    \stepcounter{Mycounter}
    {\noindent \bf \thedefinition:\ }}

\newcounter{example}[section]
\newcounter{remark}[section]
\renewcommand{\theremark}{{Remark \thesection.\arabic{remark}}}
\newcommand{\remark}{%
    \setcounter{remark}{\value{Mycounter}}
    \refstepcounter{remark}
    \stepcounter{Mycounter}
    {\noindent \bf \theremark:\ }}

\newcounter{problem}[section]
\newcounter{question}[section]
\makeatletter

\@addtoreset{equation}{section} \@addtoreset{footnote}{section}
\makeatother

\def\blacksquare{\hbox{\vrule width 5pt height 5pt depth 0pt}}
\def\endproof{\blacksquare}

\begin{document}

\begin{center}
{\Large\bf
MBM loci in families of hyperk\"ahler manifolds\\[3mm]
and centers of birational contractions\\[3mm]
}

Ekaterina Amerik\footnote{Partially supported 
by  the  Russian Academic Excellence Project '5-100'.}, 
Misha Verbitsky\footnote{Partially supported 
by  the  Russian Academic Excellence Project '5-100' and CNPq - Process 313608/2017-2.


}

\end{center}

{\small \hspace{0.15\linewidth}
\begin{minipage}[t]{0.7\linewidth}
{\bf Abstract} \\
An MBM class on a hyperk\"ahler manifold $M$
is a second cohomology class such that its orthogonal
complement in $H^2(M)$ contains a maximal dimensional
face of the boundary of the K\"ahler cone for some 
hyperk\"ahler deformation of $M$. 
An MBM curve is a rational curve in an MBM class and such that its
local deformation space has minimal possible dimension $2n-2$,
where $2n$ is the complex dimension of $M$. 
We study the MBM loci, defined as the subvarieties covered by
deformations of an MBM curve within $M$.
When $M$ is projective, MBM loci are
centers of birational contractions.
For each MBM class $z$,
we consider the Teichm\"uller space  $\Teich^{\min}_z$ of
all deformations of $M$ such that $z^{\bot}$ contains a face of the K\"ahler cone.
We prove that for all $I,J\in \Teich^{\min}_z$,
the MBM loci of $(M, I)$ and $(M,J)$ 
are homeomorphic under a homeomorphism preserving the MBM curves, 
unless possibly the Picard number of $I$ or $J$ is maximal.
\end{minipage}
}

\tableofcontents


\section{Introduction}


\subsection{Teichm\"uller spaces in hyperk\"ahler geometry}

Let $M$ be a complex manifold. Recall that 
{\bf the Teichm\"uller space} $\Teich$ of complex structures on $M$
is the quotient $\Teich:=\Comp/\Diff_0$, where $\Comp$ is the space
of complex structures (with the topology of uniform convergence
of all derivatives) and $\Diff_0$ the connected component
of the diffeomorphism group. In this paper we are interested
in the action of the mapping class group $\Diff/\Diff_0$
on $\Teich$ (see Section \ref{_Teich_Section_}). 

In our case $M$ is a compact holomorphically symplectic
manifold of K\"ahler type.\footnote{By the Calabi-Yau theorem
(\ref{_CY_Theorem_}), this is the same as a hyperk\"ahler manifold.}
We assume that $M$ has maximal holonomy (\ref{_max_holo_Definition_}; such an $M$ is also called irreducible)
and consider the Teichm\"uller space of all complex structures of
hyperk\"ahler type (Subsection \ref{_Teich_subsection_}). By a result of D. Huybrechts,
this space has finitely many connected components, and we take the one containing the 
parameter point of our given complex structure; in other words we consider the
Teichm\"uller space of all hyperk\"ahler deformations of $M$. By abuse of notation,
this space is
also denoted $\Teich$. The action of the mapping class group
$\Gamma$ (i.e. the finite-index subgroup of $\Diff/\Diff_0$ preserving our connected component) on $\Teich$ is ergodic,
and its orbits are classified using the Ratner's 
orbit classification theorem (\ref{orbits-teich}).

In \cite{_Markman:universal_}, E. Markman
has constructed the universal family
\[ u: {\cal U} \arrow\Teich.\]
The map $u$ is a smooth complex analytic
submersion with fiber $(M,I)$ at a point 
$I\in \Teich$ (throughout the paper, $(M,I)$ denotes
a manifold $M$ equipped with a complex structure $I$).
In this paper we investigate the action of the mapping
class group $\Gamma$ on this universal
fibration. We are interested in applications to the geometry of
rational curves $C\subset (M,I)$.

Fix a homology class $z\in H_2(M,\Z)$.
Let $\Gamma_z\subset \Gamma$ be the stabilizer of $z$ in
$\Gamma$, and $\Teich_z$ the Teichm\"uller space of all
complex structures $I\in \Teich$ such that
$z$ is of Hodge  type (1,1) on $(M,I)$.

Recall that the second cohomology group
 of a hyperk\"ahler
manifold with maximal holonomy is equipped
with a canonical bilinear symmetric
pairing $q$, called Bogomolov-Beauville-Fujiki (BBF) form
(\ref{_BBF_Definition_}). This form is integral but in general not unimodular,
so that it embeds $H_2(M, \Z)$ into $H^2(M,\Q)$ as an overlattice of $H^2(M,\Z)$. It is often
convenient to consider the homology classes of curves as second cohomology classes with rational
coefficients, which we do throughout the paper. In the projective case, rational curves
in the class $z$ satisfying
$q(z,z)<0$ can be contracted birationally
by the Kawamata base point free theorem
(\ref{_Kawamata_bpf_Theorem_}). It is interesting to study
the loci covered by such curves. In this paper we are
concerned with the question of invariance of these loci by
deformations.

It turns out, and is very helpful for such a study, that for any integral class $z\in H_2(M,\Z)$
with $q(z, z) <0$,  the action of the group
$\Gamma_z$ on $\Teich_z$ is also ergodic on each
connected component. Moreover we can classify, in the same way as for $\Gamma$ acting on $\Teich$,
the orbits of $\Gamma_z$-action on the space $\Teich_z^{\min}$, which is the same as $\Teich_z$ up to
inseparability issues (\ref{orbits-teichz}).

More precisely, the space $\Teich_z$
is a smooth, non-Hausdorff manifold, equipped with
a local diffeomorphism to the corresponding period
space $\Perspace_z:=\frac{SO(3, b_2-4)}{SO(1, b_2-4)\times SO(2)}$ (alternatively, this is just the orthogonal of $z$ in the usual period space $\Perspace$, seen as a subset of a quadric in the projective space $\P H^2(M, \C)$) , which becomes one-to-one if we glue together the unseparable points.
Following E. Markman \cite{_Markman:survey_}, the set of preimages of $I\in \Perspace_z$ 
(that is, the set of points unseparable from a given one $(M,I)$)  
is identified
with the set of the K\"ahler chambers in the positive cone
of $H^{1,1}(M,I)$, so that each K\"ahler chamber can be seen as the K\"ahler cone of the corresponding complex structure. For the classes relevant for us, the so-called MBM classes 
(Subsection \ref{_MBM_intro_Subsection_} and Section \ref{sectMBM}), the orthogonal complement $z^\bot$
contains one of the walls of these K\"ahler chambers
(\ref{_MBM_walls_Theorem_}). It is convenient to restrict oneself to the K\"ahler
chambers adjacent to the $z^\bot$ wall (in other words, the complex structures where $z^\bot$ is actually part of the 
boundary of the K\"ahler cone): in this way we obtain the space $\Teich_{\pm z}^{\min}\subset \Teich_z$. Note that
both spaces are non-Hausdorff even at their general points, since there are always at least two chambers adjacent to a 
given wall. We remark that $z^\bot$ is co-oriented, and pick the set of the
chambers adjacent to $z^\bot$ on the positive side.
This last space, separated at its general point, is denoted $\Teich_z^{\min}\subset \Teich_z$.  
This is precisely the space of complex structures
$I\in \Teich_z$ such that a positive multiple of $z$ is represented by an extremal rational curve: indeed, by a result of
Huybrechts and Boucksom, the K\"ahler cone is
characterized as the set of $(1,1)$-classes which are
positive on all rational  curves (\cite{_Huybrechts:basic_,_Huybrechts:erratum_}, \cite{_Boucksom-cone_}).

In this paper we are interested in subvarieties
$Z_I\subset M_I$ covered by extremal or minimal curves.
It is known since \cite{_Ran:deformations_} (see also
\cite{_AV:MBM_}) that a minimal curve in the cohomology 
class $z$ on a hyperk\"ahler manifold $(M,I)$ deforms locally together with $(M,I)$
as long as $z$ remains of type $(1,1)$. Since $\Teich_z$ is not Hausdorff, we cannot deduce from this that the curve
deforms as long as the complex structure stays in $\Teich_z\subset \Teich$. However this is obviously true modulo this
inseparability issue: for each complex structure $I\in \Teich_z$ there exists an $I'$ unseparable from $I$ and carrying
a deformation of our curve (which, however, does not have to be extremal anymore - it can, for instance, become reducible).
It turns out that we can actually be more precise: over $\Teich_z^{\min}$, the homeomorphism class (and actually stratified
diffeomorphism, and also bi-Lipschitz equivalence class)
of the locus $Z_I\subset (M,I)$ defined as above does not depend on $I$,
except possibly for the complex structures with maximal Picard number.


The proof of this is based on the following
observation. Let $E\stackrel \phi\arrow  B$ be a proper 
holomorphic  (or even real analytic) map, and
assume that $B$ is obtained as
a union of dense subsets, $B= \coprod_{\alpha \in {\goth
    I}} B_\alpha$, such that for any index $\alpha$,
all fibers of $\phi$ over $b\in B_\alpha$ are isomorphic. Then
all fibers of $\phi$ are homeomorphic, stratified
diffeomorphic and bi-Lipschitz equivalent.

This observation is based on the classical results
by Thom and Mather (the bi-Lipschitz case is
due to Parusi\'nski\footnote{\cite{_Parusinski:Lip_}, 
\cite{_Parusinski:Lip_stra_}.}). They proved that for any
proper real analytic fibration $E\stackrel \phi\arrow B$,
there exists a stratification of $B$ such that the
restriction of $\phi$ to open strata is locally
trivial in the category of topological spaces
(or in bi-Lipschitz category). Since each $B_\alpha$
in the decomposition $B= \coprod_{\alpha \in {\goth
    I}} B_\alpha$ intersects the open stratum,
this implies that all fibers of $E\stackrel \phi\arrow B$
are homeomorphic and bi-Lipschitz equivalent.

The dense subsets $B_\alpha$ are in our case provided by the
ergodicity of the mapping class group action.

We state our main results precisely in the subsection \ref{main}, after a brief digression on rational
curves in the next subsection.

\subsection{MBM loci on hyperk\"ahler manifolds}
\label{_MBM_intro_Subsection_}

Let $C\subset M$ be a rational curve on a holomorphic symplectic manifold of dimension $2n$. According to a theorem of Ran \cite{_Ran:deformations_},
the irreducible components of the deformation space of $C$ in $M$ have dimension at least $2n-2$.

\hfill

\definition\label{minimal}
A rational curve $C$ in a holomorphic
symplectic manifold $M$ is called {\bf minimal}
if every component of its deformation space has dimension $2n-2$ at $C$.

\hfill

The dimension of a maximal irreducible uniruled subvariety of $M$ can take any value between $n$ and $2n-1$. Such a subvariety is always coisotropic, and applying bend-and-break lemma one sees that there is always a minimal curve through a general point 
of such a subvariety (\cite{_AV:MBM_}, Section 4).

\hfill

The key property of a minimal curve is that such a curve deforms together with its cohomology class $[C]$. More precisely,
any small deformation of $M$ on which $[C]$ is still of type $(1,1)$, contains a deformation of $[C]$ 
(\cite{_AV:MBM_}, Corollary 4.8). Taking closures in the
universal family over $\Teich_z$ gives a submanifold of
$\Teich_z$ of maximal dimension (which does not have to
coincide with $\Teich_z$, as it is not Hausdorff) such that every complex structure in this 
submanifold carries a deformation of $C$; this curve,
however, can degenerate to a reducible curve,
and one cannot in general say much about the
cohomology classes of its components (Markman's example on
K3 surfaces is already enlightening, see
\cite{Markman-pex}, Example 5.3).

\hfill

In \cite{_AV:MBM_}, we have defined and studied the MBM classes: these are classes $z\in H^2(M,\Z)$ such that, up to monodromy and 
birational equivalence, $z^{\bot}$ contains a face of the
K\"ahler cone\footnote{By convention, our ``face'', or
  ``wall'', shall always mean a face of maximal dimension,
  which is equal to $h^{1,1}-1$.}. In other words, $z^{\bot}$ contains a face of some K\"ahler chamber (see \cite{_Markman:survey_} for the definition of the latter, but it amounts to say that those are monodromy transforms of K\"ahler cones of the birational models of $M$). It is clear that 
the Beauville-Bogomolov square $q(z)$ is then negative; on the other hand, one can characterize MBM classes as negative classes 
such that some rational multiple $\lambda z$ is represented by a rational curve on a deformation of $M$ (\cite{_AV:MBM_}, Theorem 5.11). For our purposes, it is convenient to extend the notion of MBM on the rational cohomology (or integral homology) classes in an obvious way.

We would like to study rational curves whose class is MBM and prove certain deformation invariance statements related to such curves. The notion of MBM class is defined up to a rational multiple, whereas by the above discussion we want 
to restrict to ``minimal'' rational curves. The most straightforward definition of minimality, cf. \cite{_AV:MBM_}, is probably ``of minimal degree in the
uniruled subvariety covered by its deformations''; then by bend-and-break the curve deforms within $M$ in a family of dimension
at most $2n-2$ and one deduces from this as in \cite{_Ran:deformations_} or \cite{_AV:MBM_} that it deforms outside $M$ together with its cohomology class. However the uniruled
subvariety in question might have several irreducible components. To keep track of this we have to consider the minimality with respect to those
components, hence \ref{minimal}.

Note that it is apriori possible (though we don't have any example) that the same rational curve $C$ is contained in two maximal irreducible uniruled subvarieties $Z_1$ and $Z_2$ of $M$, in such a way that the deformations of $C$ lying in $Z_1$ form a $2n-2$-parameter family whereas those 
lying in $Z_2$ need more parameters. Such a $C$ is, by our definition, not minimal, but its generic deformation in $Z_1$ is.

\hfill

\definition
{\bf An MBM curve} is a minimal curve $C$ such that its class $[C]$ is MBM.

\hfill

\definition
Let $C$ be an MBM curve on a hyperk\"ahler manifold $(M,I)$,
and $B$ an irreducible component of its deformation
space in $M$ containing the parameter point for $C$. An {\bf MBM locus} of $C$ is the union of all curves
parameterized by $B$.

\hfill

As mentioned in the beginning of this subsection, the MBM loci are coisotropic subvarieties which can have any dimension between $n$ and $2n-1$, but the family of minimal 
rational curves in an MBM locus always has $2n-2$ parameters.

\hfill

\definition\label{full-MBM} Let $z$ be an MBM class in $H^2(M,\Q)$. The {\bf full MBM locus} of $z$ is the union of all MBM curves of cohomology class proportional to $z$ and their degenerations (in other words, the union of all MBM loci for MBM curves of cohomology class proportional to $z$).

\hfill

\remark If the complex structure on $M$ is in
$\Teich_z^{\min}$, the full MBM locus has only finitely
many irreducible components and is simply the union of all
rational curves of cohomology class proportional to
$z$. We sketch the argument here and refer to
Section \ref{_MBM_loci_Section_} for details.
The reason for the finiteness is that when $(M,I)$ is projective, this
MBM locus is the exceptional set of a birational
contraction. Indeed, for any rational point in the interior of the face
of the K\"ahler cone associated with $[C]^\bot$,
the corresponding line bundle is semiample
by the Kawamata base-point-free theorem, and the corresponding
birational contraction contracts exactly the curves which have cohomology class proportional to that of $C$. The number of irreducible components
of an exceptional set of a contraction is finite. One knows that these are uniruled and by bend-and-break lemma one finds a minimal rational curve in each.

On the other hand, any $(M,I)$ admits a projective 
small deformation. By semi-continuity one deduces that the number of irreducible components
of a locus covered by rational curves of cohomology class proportional to $z$ is always finite.
Finally, there is a minimal curve in every uniruled component, and since we are in $\Teich_z^{\min}$,
the class of such an MBM curve remains proportional to $z$, so that our uniruled components are MBM loci.

\hfill

\subsection{Main results of this paper}\label{main}


Global deformations of rational curves, even minimal ones, in a given cohomology class $z\in H^2(M,\Q)$ (where 
$M$ is the underlying differentiable manifold) are rather 
difficult to understand. 
As an illustration of what we still cannot do, consider a general projective 
K3 surface $S$: by an argument
due to Bogomolov and Mumford
(\cite{MM} and more recent work \cite{BoHT}, \cite{LL}), $S$ contains (singular)
rational curves. It follows that the second
punctual Hilbert scheme of $S$ (often denoted by $S^{[2]}$) contains a rational surface (swept out by 
deformations of the diagonal in the symmetric square of a rational curve). A natural question is
whether each projective deformation of $S^{[2]}$ (the deformations of $S^{[2]}$ are often called {\bf irreducible holomorphic symplectic 
fourfolds of K3 type}) contains a rational surface. For the moment the answer seems to be still unknown.
A natural idea would be to deform the surface together with the cohomology class of a minimal rational 
curve in it, however this does not always work.
For instance, the surface may be contained in a uniruled
divisor swept out by deformations of the same minimal
curve, in such a way that only the divisor  
survives on the neighbouring manifolds.

The main point of the present paper is that it turns out to be much easier
to understand the deformations of 
subvarieties swept out by minimal rational curves with
negative Beauville-Bogomolov square. As we have mentioned before,
the Beauville-Bogomolov quadratic form induces an
embedding from $H_2(M,\Z)$ to $H^2(M,\Q)$ so that
we may view the homology classes of curves as rational
cohomology classes of type $(1,1)$. Therefore, 
it makes sense to talk of the Beauville-Bogomolov square
of a curve. By \cite{_AV:MBM_}, 
the class $z$ of such a curve is MBM, meaning that up to monodromy action and birational equivalence
its orthogonal supports a wall of the K\"ahler cone and this is also the case  on each deformation 
where $z$ remains
of type $(1,1)$. 


We restrict ourselves to the space $\Teich_z^{\min}\subset \Teich_{\pm z}^{\min}\subset \Teich_z$ described in  
the first subsection. Recall that to construct $\Teich_z^{\min}\subset \Teich_z$, we
first take the complex structures where $z$ actually contains a wall of the K\"ahler cone obtaining $\Teich_{\pm z}$, then
take the ``positive half'' of it. 
 On the space $\Teich_z^{\min}$, there is an action of the subgroup of
 the monodromy group preserving $z$, and 
it turns out, thanks to the negativity of $z$, that  
almost all orbits of this action are dense. This allows us
to make conclusions such as the uniform behaviour of subvarieties swept out by curves of class $z$
on the manifolds represented by the points of $\Teich_z^{\min}$.

Our main result is as follows.

\hfill

\theorem\label{_main_homeo_Theorem_}
Let $M$ be a hyperk\"ahler manifold of maximal holonomy and $z\in H_2(M,\Z)\subset H^2(M,\Q)$ a class of negative 
Beauville-Bogomolov square. Assume that $z$ is represented by a minimal rational curve
in some complex structure $I$ on $M$ (this means that $z$ and the curve are MBM, see \cite{_AV:MBM_}).
Let $Z=Z_I\subset (M,I)$ be the union of all rational
curves in  cohomology classes proportional to $z$ and degenerations of
such curves.
Then for all $I,I'\in \Teich_z^{\min}$ such that the Picard
number of $(M,I)$  and $(M,I')$ is not maximal (that is, not equal
to $h^{1,1}(M_I)$), there exists a homeomorphism
$h:\; (M,I) \arrow (M,I')$ identifying $Z_I$ and $Z_{I'}$.

\proof See \ref{homeo-loci}. \endproof

\hfill

\remark
In fact $h$ is more than just a 
homeomorphism: the spaces $Z_I$ and $Z_{I'}$ are naturally stratified 
by complex analytic subvarieties, and $h$ is a
diffeomorphism on 
open strata.

\hfill

\remark 
Any complex variety $X$ can be locally embedded to $\C^n$.
Consider the path metric on $X$ obtained from the usual
metric on $\C^n$. It is not hard to see that the
bi-Lipschitz class of this metric is independent
from the choice of the local embedding. The homeomorphism
$h$ constructed in \ref{homeo-loci} 
is bi-Lipschitz (Subsection \ref{_Variants_subsection_}). 

\hfill

\remark
When the manifold $(M,I)$ is projective,
the varieties $Z_I\subset (M,I)$ 
are exceptional loci of birational contractions. This observation
follows directly from the Kawamata
base point free theorem (\ref{_Kawamata_bpf_Theorem_}). Kawamata 
base point free theorem in non-algebraic
setting is unknown, but we conjecture that
$Z_I$ are centers of bimeromorphic contractions 
for non-algebraic K\"ahler deformations of $(M,I)$ as well.

\hfill

In Subsection \ref{_Variants_subsection_}
we explain the following two variants/strengthenings of
\ref{_main_homeo_Theorem_}.

\hfill

\theorem\label{_main_incidentn_Theorem_}
In assumptions of \ref{_main_homeo_Theorem_},
let $B_I$ be the Barlet space of all rational
curves of cohomology class proportional to $z$. Then
the homeomorphism $h:\; Z_I \arrow Z_{I'}$
can be chosen to send any rational curve $C\in B_I$ to some rational curve $h(C)\in B_{I'}$,
inducing a homeomorphism from $B_I$ to $B_{I'}$.

\hfill

Recall that a compact K\"ahler manifold has a so-called
MRC fibration  (\cite{_Campana:MRC_}, \cite{_KMM:MRC_})
whose fiber at a 
general point $x$ consists of all the points which can be reached from $x$ by a chain of
rational curves. In particular, considering such a fibration on a desingularization of a
component of $Z_I$ gives a rational map $Q:\; Z_I\dasharrow Q_I$.

\hfill

\theorem\label{_main_MRC_Theorem_}
In assumptions of \ref{_main_homeo_Theorem_},
consider the MRC fibrations $Q:\; Z_I\dasharrow Q_I$, $Q':\; Z_{I'}\dasharrow Q_{I'}$
Then $h:\; Z_I \arrow Z_{I'}$ induces a biholomorphism between open dense
subsets of the fibers
of $Q$ and $Q'$.

\hfill

\remark
Notice that a homeomorphism between normal
complex analytic spaces which is holomorphic
in a dense open set is holomorphic everywhere.
This result follows from a version of Riemann
removable singularities theorem, see e. g. 
\cite[Theorem 1.10.3]{_Magnusson:cycle_}.

\hfill

When the complex structure we consider is in
$\Teich_z^{\min}$, the set $Z$ defined as above is the
full MBM locus of $z$. We also have a similar statement
for MBM loci of curves, which is a straightforward
consequence of \ref{_main_homeo_Theorem_} and
\ref{_main_incidentn_Theorem_}.

\hfill

\theorem
Let $z$ be an MBM class on a hyperk\"ahler manifold, $C\subset (M,I)$ an MBM curve in this class
and $Z_C$ its MBM locus.
Then $C$ is deformed to an MBM curve $C_J \subset (M,J)$
 for all $J\in \Teich_z^\min$, and 
the corresponding MBM locus $Z_{C_J}$ 
is homeomorphic to the MBM locus $Z_C$, except possibly if the Picard number of $(M,I)$ or $(M,J)$ is maximal.
This homeomorphism  
can be chosen in such a way that all MBM
deformations of $C$ in $Z_{C}$ are mapped to 
MBM deformations of $C_J$ in $Z_{C_J}$.

\hfill

\remark
It is conceivable that two homology classes
$kz$ and $lz$, where $k$ and $l$ are positive rational numbers, would both be represented by
an MBM curve in the same manifold. Of course, the corresponding MBM
loci would be different. So far no such example is known.

\hfill

The main theorem shall be proved in 
Subsection \ref{_Proofs_Subsection_},
and the Subsection \ref{_Variants_subsection_} 
is devoted to the variants.

\hfill

\remark One cannot affirm that the same statements hold along the whole of $\Teich_z$, and this is
false already for K3 surfaces. Indeed a $(-2)$-curve on a K3 surface $X$ can become reducible on a suitable 
deformation $X'$. What we do affirm is that in $\Teich_z$ there is another point, nonseparable from the one
corresponding to $X'$, such that on the corresponding K3 surface $X''$ our curve remains irreducible.
In this two-dimensional case, this easily follows from the description of the decomposition into the K\"ahler chambers in
\cite{_Markman:survey_}; \ref{_main_homeo_Theorem_} allows us to go further in the
higher-dimensional case, unless possibly if $X'$ is of maximal Picard rank.









\hfill









\section{Hyperk\"ahler manifolds}


\subsection{Hyperk\"ahler and holomorphically symplectic manifolds}

Here we remind basic results and definitions of
hyperk\"ahler and holomorphically symplectic
geometry. Please see \cite{_Besse:Einst_Manifo_} 
and \cite{_Beauville_} for more details and reference.

\hfill

\definition
A {\bf hyperk\"ahler structure} on a manifold $M$
is a Riemannian structure $g$ and a triple of complex
structures $I,J,K$, satisfying quaternionic relations
$I\circ J = - J \circ I =K$, such that $g$ is K\"ahler
for $I,J,K$.

\hfill

\remark 
A hyperk\"ahler (i.e. the one carrying a hyperk\"ahler structure) manifold  has three symplectic forms
\[ \omega_I:=  g(I\cdot, \cdot), \omega_J:=  g(J\cdot, \cdot),
\omega_K:=  g(K\cdot, \cdot).
\]

\hfill







\definition
A holomorphically symplectic manifold 
is a complex manifold equipped with nowhere degenerate holomorphic
$(2,0)$-form.

\hfill

\remark 
Hyperk\"ahler manifolds are holomorphically symplectic.
Indeed, $\Omega:=\omega_J+\1\omega_K$ is a holomorphic symplectic
form on $(M,I)$.

\hfill

\theorem \label{_CY_Theorem_}
(Calabi-Yau, \cite{_Yau:Calabi-Yau_}; see \cite{_Besse:Einst_Manifo_}) 
A compact, K\"ahler, holomorphically symplectic manifold
admits a unique hyperk\"ahler metric in any K\"ahler class.

\hfill

For the rest of this paper, we call a compact K\"ahler complex manifold
{\bf hyperk\"ahler} if it is holomorphically symplectic.

\hfill

\definition \label{_max_holo_Definition_}
Such a manifold $M$ is moreover called
{\bf of maximal holonomy}, or {\bf simple}, or 
{\bf IHS} (irreducible holomorphically symplectic) if $\pi_1(M)=0$, $H^{2,0}(M)=\C$.

\hfill

\theorem
(Bogomolov's decomposition, \cite{_Bogomolov:decompo_}, \cite{_Beauville_})\\ Any 
hyperk\"ahler manifold admits a finite covering
which is a product of a torus and several 
hyperk\"ahler manifolds of maximal holonomy.

\hfill

\remark 
Further on, all hyperk\"ahler manifolds
are tacitly assumed to be of maximal holonomy.

\subsection{Bogomolov-Beauville-Fujiki form}

\theorem (Fujiki)
Let $\eta\in H^2(M)$, and $\dim M=2n$, where $M$ is
hyperk\"ahler. Then $\int_M \eta^{2n}=cq(\eta,\eta)^n$,
for some primitive integer quadratic form $q$ on $H^2(M,\Z)$,
and $c>0$ a rational number.

\proof \cite{_Fujiki:HK_}. \endproof

\hfill

\definition  \label{_BBF_Definition_}
This form is called
{\bf Bogomolov-Beauville-Fujiki form}. It is defined
by the Fujiki's relation uniquely, up to a sign. The sign is determined
from the following formula (Bogomolov, Beauville; see
\cite{_Beauville_}):
\begin{align*}  \lambda q(\eta,\eta) &=
   \int_X \eta\wedge\eta  \wedge \Omega^{n-1}
   \wedge \bar \Omega^{n-1} -\\
 &-\frac {n-1}{n}\left(\int_X \eta \wedge \Omega^{n-1}\wedge \bar
   \Omega^{n}\right) \left(\int_X \eta \wedge \Omega^{n}\wedge \bar \Omega^{n-1}\right)
\end{align*}
where $\Omega$ is the holomorphic symplectic form, and 
$\lambda>0$.

\hfill

\remark The form $q$ has signature $(3,b_2-3)$.
It is negative definite on primitive forms, and positive
definite on $\langle \Omega, \bar \Omega, \omega\rangle$,
 where $\omega$ is a K\"ahler form. 


\section{Teichm\"uller spaces and global Torelli theorem}
\label{_Teich_Section_}

In this section, we recall the global Torelli theorem
for hyperk\"ahler manifolds, and state some of its applications.
We follow \cite{_Verbitsky:ergodic_} and \cite{_V:Torelli_}.

\subsection{Teichm\"uller spaces and the mapping class
  group}
\label{_Teich_subsection_}

\definition 
Let $M$ be a compact complex manifold, and 
$\Diff_0(M)$ a connected component of its diffeomorphism group
({\bf the group of isotopies}). Denote by $\Comp$
the space of complex structures on $M$, and let
$\Teich:=\Comp/\Diff_0(M)$. We call 
it {\bf the Teichm\"uller space} of complex structures on $M$

\hfill

\theorem (Bogomolov-Tian-Todorov)\\
Suppose that $M$ is a Calabi-Yau manifold. Then 
$\Teich$ is a complex manifold, possibly non-Hausdorff.

\hfill

\proof This statement is essentially contained in 
\cite{_Bogomolov:defo_}; see \cite{_Catanese:moduli_} for more details. 
\endproof

\hfill

Working more specifically with hyperk\"ahler manifolds,
one usually takes for $\Teich$ the Teichm\"uller space of all complex 
structures of hyperk\"ahler type. It is an open subset
in the Teichm\"uller space of all complex structures 
by Kodaira-Spencer K\"ahler stability theorem
\cite{_Kod-Spen-AnnMath-1960_}.

\hfill

\definition 
Let $\Diff(M)$ be the group of 
diffeomorphisms of $M$. We call $\Gamma=\Diff(M)/\Diff_0(M)$ {\bf the
mapping class group}. 

\hfill

\remark The quotient $\Teich/\Gamma$ 
is identified with the set of equivalence classes of complex structures.

\hfill

If $M$ is IHS, the space $\Teich$ has finitely many
connected components by a result of Huybrechts
(\cite{_Huybrechts:finiteness_}). We
consider the subgroup $\Gamma_0$
of the mapping class group which preserves the one
containing the parameter point for our chosen complex
structure.

\hfill

\theorem (\cite{_V:Torelli_})\\
Let $M$ be a simple hyperk\"ahler manifold,
and $\Gamma_0$ as above.  Then 
\begin{description}
\item[(i)]  The image of $\Gamma_0$ in $\Aut{H^2(M,\Z)}$ is a finite index 
subgroup of the orthogonal lattice $O(H^2(M, \Z), q)$.
\item[(ii)] The map $\Gamma_0\arrow O(H^2(M, \Z), q)$
has finite kernel.
\end{description}


\hfill

\definition We call the image of $\Gamma_0$ in
$\Aut{H^2(M,\Z)}$ the {\bf monodromy group},
denoted by $\Mon(M)$.

\hfill

\remark From now on, by abuse of notation and since we are interested in deformations of $(M,I)$, 
we denote by $\Teich$ the connected component of the Teichm\"uller space 
containing the parameter point for our given complex structure, and the mapping class group means the finite index
subgroup of $\Diff(M)/\Diff_0(M)$ preserving this component.

\subsection{The period map}



\remark
For any $J\in \Teich$,
$(M,J)$ is also a simple hyperk\"ahler manifold, hence
$H^{2,0}(M,J)$ is one-dimensional. 
 
\hfill

\definition 
Let  $\Per:\; \Teich \arrow {\Bbb P}H^2(M, \C)$
map $J$ to a line $H^{2,0}(M,J)\in {\Bbb P}H^2(M, \C)$.
The map $\Per:\; \Teich \arrow {\Bbb P}H^2(M, \C)$ is 
called {\bf the period map}.

\hfill

\remark 
From the properties of BBF form, it follows that
$\Per$ maps $\Teich$ into an open subset of a 
quadric, defined by
\[
\Perspace:= \{l\in {\Bbb P}H^2(M, \C)\ \ | \ \  q(l,l)=0, q(l, \bar l) >0\}.
\]
It is called {\bf the period space} of $M$.

\hfill

\remark 
One has
\[ \Perspace=\frac{SO(b_2-3,3)}{SO(2) \times SO(b_2-3,1)}=\Gr_{++}(H^2(M,\R).
\]
Indeed, the group $SO(H^2(M,\R),q)=SO(b_2-3,3)$ acts transitively on
$\Perspace$, and $SO(2) \times SO(b_2-3,1)$ is a stabilizer of a point.

\hfill

\theorem (Bogomolov)\\
For any hyperk\"ahler manifold, 
period map is locally a diffeomorphism.

\proof \cite{_Bogomolov:defo_}. \endproof

\subsection{Birational Teichm\"uller moduli space}

\definition
Let $M$ be a topological space. We say that $x, y \in M$
are {\bf non-separable} (denoted by $x\sim y$)
if for any open sets $V\ni x, U\ni y$, $U \cap V\neq \emptyset$.

\hfill




By a result of Huybrechts (\cite{_Huybrechts:basic_}), 
any two non-separable points $I,I'$ in the Teichm\"uller space
correspond to birational complex manifolds  $(M,I)$ and $(M,I')$. The birational map in question,
though, might be biregular: indeed the Teichm\"uller space is non-separated even for K3 surfaces.
The precise description of non-separable points of $\Teich$ can be found in \cite{_Markman:survey_}
and is as follows. Consider the positive cone $\Pos(M, I)$
which is one of the two connected components of the set of
positive vectors
\[ \{x\in H^{1,1}(M, I)|q(x,x)>0\}\] 
containing the K\"ahler cone. By a result of Huybrechts and Boucksom,
the K\"ahler classes are those elements of $\Pos(M, I)$
which are positive on all rational curves (\cite{_Huybrechts:basic_,_Huybrechts:erratum_},  \cite{_Boucksom-cone_}).
It turns out that $\Pos(M, I)$ is decomposed into chambers which are K\"ahler cones of all hyperk\"ahler 
birational models of $(M,I)$ and their transforms by monodromy. The points of $\Teich$ nonseparable from $I$ 
correspond to the chambers of this  decomposition  of $\Pos(M, I)$. We shall return to this in more detail 
in Section \ref{sectMBM}.

In particular if there is no rational curve on $(M, I)$, then the K\"ahler cone of $(M, I)$ is equal to
the positive cone and $I$ is a separated point of $\Teich$. Note
that a very general hyperk\"ahler manifold has no curves at all; the ones which
contain rational curves belong to a countable union of divisors in $\Teich$. Therefore $\Teich$
is separated ``almost everywhere''.

\hfill

\definition
The space $\Teich_b:= \Teich/\sim$ is called the
{\bf birational Teichm\"uller space} of $M$, or the {\bf Hausdorff reduction} of $\Teich$.

\hfill

\theorem (Torelli theorem for hyperk\"ahler manifolds, \cite{_V:Torelli_})\\
The period map 
$\Teich_b\stackrel \Per \arrow \Perspace$ is a diffeomorphism,
for each connected component of $\Teich_b$.

\hfill

\definition Let $z$ be a class of negative square in $H^2(M,\Z)$. We call $\Teich_z$ the part of $\Teich$
consisting of all complex structures on $M$ where $z$  of type $(1,1)$.

\hfill

The following proposition is well-known (see e. g. \cite{_AV:MBM_}).

\hfill

\proposition $\Teich_z=\Per^{-1}(z^{\bot})$, 
where $z^{\bot}$ is the set of points corresponding to
lines orthogonal to $z$ in $\Perspace\subset \P H^2(M,\C)$.

\hfill

On $\Teich_z$, we have a natural action of the stabilizer
of $z$ in $\Gamma$, denoted by $\Gamma_z\subset \Gamma$.

\subsection{Ergodicity of the mapping class group action}

\definition
Let $M$ be a complex manifold, $\Teich$ its Techm\"uller
space, and $\Gamma$ the mapping class group acting on $\Teich$.
{\bf An ergodic complex structure} is a complex
structure with dense $\Gamma$-orbit.

\hfill

This term comes from the following definition and facts.

\hfill

\definition
Let $(M,\mu)$ be a space with measure,
and $G$ a group acting on $M$ preserving the measure.
This action is {\bf ergodic} if all
$G$-invariant measurable subsets $M'\subset M$
satisfy $\mu(M')=0$ or $\mu(M\backslash M')=0$.

\hfill

The following claim is well known.

\hfill

\claim
Let $M$ be a manifold, $\mu$ a Lebesgue measure, and
$G$ a group acting on $M$ ergodically. Then the 
set of non-dense orbits has measure 0.

\hfill

{\bf Proof:}
Consider a non-empty open subset $U\subset M$. 
Then $\mu(U)>0$, hence $M':=G\cdot U$ satisfies 
$\mu(M\backslash M')=0$. For any orbit $G\cdot x$
not intersecting $U$, $x\in M\backslash M'$.
Therefore, the set $Z_U$ of such orbits has measure 0. 

{\bf Step 2:} Choose a countable base
$\{U_i\}$ of topology on $M$. Then the set of 
points in dense orbits is $M \backslash \bigcup_i Z_{U_i}$.
\endproof

\hfill

\definition
{\bf A lattice} in a Lie group is a discrete
subgroup $\Gamma\subset G$ such that $G/\Gamma$ has finite
volume with respect to Haar measure.

\hfill

\theorem (Calvin C. Moore, \cite{_Moore:ergodi_})
Let $\Gamma$ be a lattice in a non-compact 
simple Lie group $G$ with finite center, and $H\subset G$ a 
non-compact subgroup. Then the left action of $\Gamma$
on $G/H$ is ergodic.
\endproof

\hfill

\theorem Let ${\Perspace}$ be a component of 
the birational Teichm\"uller space identified with the period domain, and
$\Gamma$ its monodromy group. Let $\Perspace_e$ be
a set of all points $L\in \Perspace$
such that the orbit $\Gamma\cdot L$ is dense. Then 
$Z:=\Perspace\backslash \Perspace_e$ has measure 0.

\hfill

\pstep
 Let $G=SO(3, b_2-3)$, $H=SO(2) \times SO(1, b_2-3)$.
Then $\Gamma$-action on $G/H$ is ergodic,
by Moore's theorem. 

{\bf Step 2:}
Ergodic orbits are dense, because the union 
of all non-ergodic orbits has measure 0.
\endproof

\hfill

In \cite{_Verbitsky:ergodic_} and  \cite{_Verbitsky:erratum_}, a more precise result has been established using
Ratner theory. For reader's convenience we recall the idea of proof as well.

\hfill

\theorem\label{orbits-per} 
Assume $b_2(M)\geq 5$.
Let $\Per$ and $\Gamma$ be as above, then there are three types of $\Gamma$-orbits on
$\Perspace$:

1) closed orbits which are orbits of those period planes 
\[ V_I:=\Re(H^{2,0}(M,I)\oplus H^{0,2}(M,I))\in \Perspace
\] which are rational 
in $H^2(M,\R)$
(equivalently, the corresponding complex structures $I$ are of maximal Picard number, as $H^{1,1}(M,I)$ is then also 
a rational 
subspace);

2) dense orbits which are orbits of period planes 
$V_I$ containing no
non-zero rational vectors;

3) ``intermediate orbits'':
orbits of period planes containing a single rational vector $v$.
The orbit closure then consists of all period planes
containing $v$.

\hfill

{\bf Idea of proof:} Let $G=SO^+(3, b_2-3)$ and $H=SO^+(1,b_2-3)$, so that $G/H$ is fibered
in $SO(2)$ over $\Perspace$. The group $H$ is generated by unipotents (this is why this time
we take $H=SO^+(1,b_2-3)$ rather than $H=SO(2)\times SO^+(1,b_2-3)$, so that Ratner theory
applies to the action of $H$ on $\Gamma \backslash G$,
where $\Gamma$ remains the same as above). 
Ratner
theory describes orbit closures of this action: $\overline
{xH}$ is again an orbit under a closed 
intermediate subgroup $S$, also generated by unipotents
and in which $\Gamma\cap S$ is a lattice.
From the study of Lie group structure on $G$ we derive that the subgroup must be either $H$ 
itself (the orbit is closed), or the whole of $G$ (the orbit is dense), or the stabilizer of an 
extra vector $\cong SO^+(2, b_2-3)$ (the third case). One concludes passing in an obvious way (via the double quotient)
from an $H$-action on $\Gamma \backslash G$ to a $\Gamma$-action on $G/H$. \endproof

\hfill

A useful point for us is the following observation from \cite{_Verbitsky:erratum_}.

\hfill

\proposition \label{_interme_orbit_not_anal_Proposition_}
In the last case, the orbit closure is a fixed point set
of an antiholomorphic involution, in particular, it is not
contained in any complex submanifold nor contains any
positive-dimensional complex submanifold.

\hfill

For our present purposes we need the following variant of \ref{orbits-per}.

\hfill

\theorem\label{orbits-perz} Assume $b_2(M)>5$. Let $z\in H^2(M,\Z)$ be a negative class and $\Perspace_z=z^{\bot}\subset \Perspace$ be
the locus of period points of complex structures where $z$ is of type $(1,1)$. Let $\Gamma_z$ be the subgroup of $\Gamma$ fixing $z$. Then the same conclusion as in \ref{orbits-per} holds, namely
an orbit of $\Gamma$ is either closed and consists of points with maximal Picard number, 
or dense, or the orbit closure consists of period planes containing a rational vector and in this last case it is not contained in any complex subvariety nor contains one.

\hfill

\proof It is exactly the same as the proof of \ref{orbits-per} once one interprets $\Perspace_z$ as the Grassmannian of positive 2-planes in $V=z^{\bot}\subset H^2(M,\R)$, takes 
$G=SO^+(V)\cong SO^+(3, b_2-4)$, $H=SO^+(1,b_2-4)$ and replaces $\Gamma$ by $\Gamma_z$. \endproof 





\section{MBM curves and the K\"ahler cone}
\label{sectMBM}


The notion of MBM classes was introduced in \cite{_AV:MBM_}
and studied further in \cite{_AV:Mor_Kaw_}. We recall the setting and
some results and definitions from these papers.

\hfill

\definition
The BBF form on $H^{1,1}(M,\R)$ has signature $(1, b_2-3)$.
This means that the set $\{\eta \in H^{1,1}(M,\R)\ \ |\ \ (\eta,\eta)>0\}$
has two connected components. The component which contains the
K\"ahler cone $\Kah(M)$ is called {\bf the positive cone},
denoted $\Pos(M)$.

\hfill

The starting point is the following theorem.

\hfill

\theorem
(Huybrechts \cite{_Huybrechts:basic_,_Huybrechts:erratum_}, Boucksom \cite{_Boucksom-cone_})\\
The K\"ahler cone of $M$ is the set of all $\eta\in \Pos(M)$
such that $(\eta, C)>0$ for all rational curves $C$.

\hfill

 Remark that it is sufficient to consider the curves of negative square (as 
only these have orthogonals passing through the interiour of the positive cone) and moreover
extremal, i.e. such that their cohomology class cannot be decomposed as a sum of classes of other
curves. An extremal curve is minimal in the sense of our \ref{minimal}, though apriori the converse needs 
not be true.









\hfill

\remark Here, using the BBF form, we identify
$H_2(M,\Q)$ and $H^2(M,\Q)$ for any hyperk\"ahler manifold $M$. 
This allows us to talk of the BBF square of a homology
class of a curve.

\hfill



The K\"ahler cone is thus locally polyhedral in the interior of the positive cone (with
some round pieces in the boundary), and its faces are supported on the
orthogonal complements to the extremal curves.

The notion of an extremal curve is however not adapted to the defor\-mation\--in\-variant context. 
In order to put the theory in this
context we have defined the MBM (monodromy birationally minimal) classes in \cite{_AV:MBM_}.
Here we recall several equivalent definitions (we refer to Section 5 of \cite{_AV:MBM_} for proofs of equivalences).

\hfill

\definition A negative class $z$ in the image of
$H_2(M,\Z)$ in $H^2(M,\Q)$ is called MBM if $\Teich_z$
contains no twistor curves.
This is equivalent to saying that a rational multiple of $z$ is represented by a curve in some complex 
structure where the Picard group is generated by $z$ over the rationals, and also to saying that
in some complex structure $X=(M,I)$ where $z$ is of type
$(1,1)$, the orthogonal complement $\gamma(z)^{\bot}$  contains a face of the
K\"ahler cone of a birational model $X'=(M,I')$ of $X$ (whence the terminology). Moreover in these two
equivalent definitions,
``some'' may be replaced by ''all'' without changing the content.

\hfill

\remark
Let $M$ be a hyperk\"ahler manifold. We call
{\bf a codimension 1 face}, or just 
{\bf a face} of the K\"ahler cone if there is no risk of confusion, a subset of its boundary with nonempty interior
obtained as the intersection of this boundary and a hyperplane in $H^{1,1}(M, \R)$ -
in other words, for us a face is always of maximal dimension unless otherwise specified.

\hfill

\theorem \label{_MBM_walls_Theorem_} (\cite{_AV:MBM_}, Theorem 6.2)
The K\"ahler cone is a connected component of the
complement, in $\Pos(M)$, of the union of 
hyperplanes $z^{\bot}$ where $z$ ranges over MBM classes of type $(1,1)$.

\hfill

\definition (cf. \cite{_Markman:survey_}) 
{\bf The K\"ahler chambers} are other connected components of this 
complement.

\hfill

Moreover we have the following connection between the K\"ahler chambers and the inseparable points of the
Teichm\"uller space (note that the decomposition of $\Pos$ into the K\"ahler chambers is an invariant
of a period point rather than of the complex structure itself, since it is determined by the position
of $H^{1,1}$ in $H^2(M,\R)$):

\hfill

\theorem (\cite{_Markman:survey_}), theorem 5.16) The points of a fiber of $\Perspace$ over a period point 
are in bijective 
correspondence with the K\"ahler chambers of the decomposition of the positive cone of the corresponding 
Hodge structure.

\hfill

\definition \label{_Teich_min_Definition_}
The space $\Teich^{\min}_{\pm z}\subset \Teich_z$ is obtained by removing the complex structures 
where $z$ does not support a wall of the K\"ahler cone. 

\hfill

In other words, at a general point of $\Teich_z$, where the Picard group is generated by $z$ over the rationals,
  $\Teich^{\min}_{\pm z}$ coincides with $\Teich_z$, whereas at special points of $\Teich_z$ where we have other MBM 
classes as well, we remove those complex structures where e.g. $z$ becomes a sum of two effective classes,
and rational curves representing $z$ thus cease to be extremal.

\hfill

Notice that the space $\Teich^{\min}_{\pm z}$ is not separated even at its general point, since $z^{\bot}$ divides the positive cone
in at least two chambers. In order to avoid working with such generically non-separated spaces we divide $\Teich^{\min}_{\pm z}$
in two halves: 

\hfill

\definition The space $\Teich^{\min}_{z}$ is the part of $\Teich^{\min}_{\pm z}$ where $z$ has non-negative intersection with K\"ahler classes (that is, $z$ is pseudo-effective).

\hfill

Now at a general point $\Teich^{\min}_z$ coincides with $\Perspace_z$ (but at special points it is still non-separated).


\section{Loci of MBM curves in families}\label{locifam}
\label{_MBM_loci_Section_}

\proposition Let $z$ be an MBM class in some complex structure $I\in \Teich_z^{min}$. The full MBM locus $Z$
is the union of all
rational curves $C$ such that $[C]$ is proportional to $z$.

\hfill

\proof We have defined the full MBM locus as the union of subvarieties swept out by minimal rational
curves of cohomology class proportional to $z$, so clearly the full MBM locus is included in the union of all 
rational curves of cohomology class proportional to $z$. On the other hand, take any component of the latter. By bend-and-break one can find a minimal rational curve through a general point of this component 
(see for example \cite{_AV:MBM_}, theorem 4.4, corollary 4.6), so this is also a component of $Z$.

\hfill

These loci are interesting since, at least in the projective case, these are centers of elementary birational 
contractions (Mori contractions). Indeed, recall the following partial case of Kawamata base-point-freeness
theorem.

\hfill

\theorem \label{_Kawamata_bpf_Theorem_}
(Kawamata bpf theorem, \cite{_Kawamata:Pluricanonical_})\\
Let $L$ be a nef line bundle on a projective manifold $M$ such that $L^{\otimes a}\otimes{\cal O}(- K_M)$ is big
for some $a$.
Then $L$ is semiample.

\hfill

Here a line bundle $L$ is said to be {\bf nef} if $c_1(L)$ is in the closure of the K\"ahler cone, and {\bf big}
if the dimension of the space of global sections of its tensor powers has maximal possible growth. For the nef
line bundles this last condition is equivalent to 
$c_1(L)^{\dim M}>0$ (that is, the maximal self-intersection 
number of $L$ being positive). A {\bf semiample} line bundle is 
a line bundle $L$ such that $L^{\otimes n}$ is base point free
for some $n$; then for $n$ big enough the linear system of sections of $L$ defines a projective morphism with
connected fibers $\phi: M\to M_0$. The bigness of $L$ in
fact implies that $\phi$ is birational. Clearly, for a
curve $C$,  $\phi(C)$ is a point if and only if $L\cdot C=0$.

\hfill

If $M$ is a holomorphic symplectic manifold, the canonical divisor $K_M$ is zero and any big and nef line
bundle is semiample. Let $z$ be a class in $H^2(M,\Z)$ such that $z^{\bot}$ contains a face of the K\"ahler 
cone. If $M$ is projective, there is an integral point in the interior of this face, and this point 
is the Chern class of a nef and big line bundle $L$. By Kawamata base point freeness $L$ is semiample
and the morphism associated with the space of sections of $L^{\otimes n}$ contracts exactly the curves
with cohomology class proportional to $z$. It is well-known that the exceptional set of a birational contraction on a hyperk\"ahler manifold is covered by rational curves (one can deduce this
for instance from \cite{Kawamata}, Theorem 1). It follows that the exceptional set of $\phi$ is exactly the
locus $Z$. 

\hfill

\corollary $Z$ has finitely many irreducible components.

\hfill

\proof Minimal rational curves survive on the small deformations of $M$ provided that their cohomology class stays of type $(1,1)$ (see \cite{_AV:MBM_}, corollary 4.8), and their loci can only collide over closed subsets of the parameter space.
 So if $Z$ has infinitely many 
components, so does the ``neighbouring'' full MBM locus $Z_{I'}$, $I'\in \Teich_z^{\min}$. But one
can always find an $I'$ close to $I$ such that $(M, I')$ is projective,
and in this case it is an exceptional set of a holomorphic
birational contraction.

\hfill

The purpose of this section is to prove the following result ( \ref{_main_homeo_Theorem_} from the Introduction).

\hfill

\theorem\label{homeo-loci} 
Let $M$ be a simple hyperk\"ahler manifold and $z$ an MBM
class. 
For any complex structure $I\in \Teich_z^{\min}$, let $M_I$
denote the complex manifold $M$ equipped with $I$
and $Z_I\subset M_I$ the subvariety covered by rational curves of cohomology class proportional to $z$.
Then for all $I\in \Teich_z^{\min}$ such that the Picard number of $M_I$ is not maximal, the subvarieties $Z_I$ 
are homeomorphic.

\hfill

The proof uses two main ingredients:
ergodicity of monodromy action on $\Teich_z^{\min}$ and Whitney 
stratification.

\subsection{Mapping class group action on $\Teich_z^{\min}$}

The group $\Gamma_z\subset \Gamma$ obviously acts on $\Teich_z^{\min}$. Indeed the action of any 
$\gamma\in \Gamma$ is just the transport of the complex structure; if $z^{\bot}$ contains a wall of the
K\"ahler cone in a complex structure $I$, then so does $\gamma z$ in the complex structure $\gamma I$.
Notice that the same remark applies to rational curves: $\gamma C$ is a rational curve in the structure
$\gamma I$ and the minimality is preserved. So the locus $Z\subset X=(M,I)$
is sent by an element of $\Gamma_z$  to $Z_{\gamma I} \subset X'=(M, \gamma I)$. 

It turns out that the results on the mapping class group action on $\Perspace$ ``lift'' to those on the action on 
$\Teich$, but if we want to work on a subspace where $z$ remains of type $(1,1)$ this has to be $\Teich_z^{\min}$
rather than $\Teich$.

The following theorem from \cite{_Verbitsky:ergodic_}, \cite{_Verbitsky:erratum_} strengthens \ref{orbits-per}.

\hfill

\theorem\label{orbits-teich} 
Assume $b_2(M)\geq 5$. Let $\Gamma$ denote the mapping class group.
Then there are three types of $\Gamma$-orbits on $\Teich$: 
closed (where the period planes are rational,
thus the complex structures have maximal Picard number), dense (where the period planes contain no rational
vectors) and such that the closure is formed by points whose period planes contain a fixed rational vector $v$.
In the last case, the orbit closure $C_v$ is totally real, so that no neighbourhood of a point $c\in C_v$ in 
$C_v$ is contained in a proper complex subvariety of $\Teich$.

\hfill

The proof proceeds by establishing that the period map commutes with taking orbit closures, in the following 
way. Introduce the space $\Teich_K$ which consists of pairs $(I,\omega)$ where $I\in \Teich$ and $\omega\in \Kah(I)$ 
is of square 1. Calabi-Yau theorem
(\ref{_CY_Theorem_}) immediately implies that 
this is the Teichm\"uller space of pairs
(complex structure, hyperk\"ahler metric compatible with it).
As shown in \cite{_AV:Teichmuller_}, the period map
is injective on the space of hyperk\"ahler metrics;
therefore, it is injective on $\Teich_K$. In other words, $\Teich_K$ 
is naturally embedded in $\Per_K$, the homogeneous manifold
of all pairs consisting of a period point $I\in
\Perspace$ and an
element $\omega$ of square one in its positive cone (which
indeed depends only on the period point, not on the
complex structure itself). The latter is a homogeneous
space, so we can try to apply Ratner theory to prove the
following result, which clearly implies what we 
need: for any $I$, the closure of the $\Gamma$-orbit of $(I, \Kah(I))\subset \Teich_K \subset \Per_K$ contains the orbit of $(\Per I, \Pos(I))$ (here by an orbit of the subset we mean the union of its translates).
Now one can construct orbits of one-parameter subgroups which are entirely contained in $(I, \Kah(I))$ and such that the closure
of their projection to $\Per_K/\Gamma$ contains the projection of the positive cone. Indeed, one deduces from the non-maximality of the Picard number that $\Kah(I)$ has a ``round part'', for instance in the following sense: in the
intersection of $\Kah(I)$ with a general 3-dimensional
subspace $W$ in $H^{1,1}(I)$, of signature $(1,2)$. This
is used to find many horocycles in $\Kah(I)$
tangent to the round part of the boundary. 
The horocycle is an orbit of a one-parameter unipotent subgroup.
Applying Ratner theory to a sufficiently general horocycle
of this type, one sees that the closure of its image in
$\Per_K/\Gamma$ contains an entire $SO(H^{1,1}(I))$-orbit
(\cite{_Verbitsky:erratum_}, Proposition
3.5), which is the positive cone $\Pos(I)$.

\hfill

The analogue of \ref{orbits-teich} in our setting is as follows.

\hfill

\theorem\label{orbits-teichz}
Assume $b_2(M)>5$ and let $z\in H^2(M,\Z)$ be an MBM class and $\Gamma_z$
the subgroup of the mapping class group consisting of all elements whose action on the second cohomology fixes $z$. Then $\Gamma_z$
acts on $\Teich_z^{\min}$ ergodically,
and there are the same three types of orbits
of this action as in \ref{orbits-teich}.

\hfill

\proof It proceeds along the same lines. We introduce the
spaces $\Per_{K,z}$ consisting of pairs 
\[
\{(\Per(I),
\omega\in \Pos(I)\cap z^{\bot}), q(\omega, \omega)=1\}
\]  
and $\Teich_{K,z}$ consisting of pairs $(I\in
\Teich_z^{\min}, \omega)\in \Per_{K,z}$ where $\omega$ 
belongs to the wall of $\Kah(I)$ given by $z^{\bot}$.
We denote such a wall by  $\Kah(I)_z$, though of course 
its elements are not K\"ahler forms on $I$, but rather
semipositive limits of those. Since the complex
structures in $\Teich_z^{\min}$ which 
have the same period point are in one-to-one
correspondence with the walls of the K\"ahler chambers in
which the other MBM classes partition $z^{\bot}$, $\Teich_{K,z}$
again embeds naturally in $\Per_{K,z}$. We fix a complex
structure $I$ with non-maximal Picard number. We need to
prove that the closure of the $\Gamma_z$-orbit of 
the subset $(I,\Kah(I)_z)$ contains the orbit of
$(\Per(I), \Pos(I)\cap z^{\bot})$. This is done exactly in
the same way as in \ref{orbits-teich}. We take a general
three-dimensional subspace $W$ in $z^{\bot}$, the
intersection of $W$ with our wall
$\Kah(I)_z$ contains horocycles, and we deduce from Ratner
orbit closure theorem and Proposition 3.5 of
\cite{_Verbitsky:erratum_} that the closure of the projection
of such a horocycle to $\Per_{K,z}/\Gamma_z$ is large, containing an $SO(H^{1,1}(I)\cap z^{\bot}$-orbit, which is the
projection of  $\Pos(I)\cap z^{\bot}$.
\endproof

\subsection{Stratification}

Consider the universal family ${\cal X}$ over
$\Teich_z^{\min}$ (\cite{_Markman:universal_}).
We are interested in the family ${\cal Z}\subset {\cal
  X}$ with the fiber over $I\in \Teich_z^{\min}$ obtained
as the full MBM locus of $z$ on the complex manifold
$X=(M,I)$. This family can be constructed, for instance,
by taking the image of the evaluation map for all
components of the relative Barlet space corresponding to
cohomology classes proportional to $z$ and dominating 
$\Teich_z^{min}$. As $\Teich_z^{\min}$ is not Hausdorff,
we shall, whenever necessary, restrict both families to a 
small neighbourhood $U$ of some point $x$, or to a small compact 
$K$ within $U$, and denote by ${\cal X}_U$, ${\cal Z}_U$ 
the restrictions of these families. 

It is well-known that an analytic subset $W$ of a complex manifold $Y$ admits a ``nice'' stratification
({\bf Whitney stratification} or {\bf Thom-Mather
  stratification}; see 
\cite{_Mather:stability_}). Recall also the following {\bf first isotopy
lemma} by Thom (we refer to \cite{_Mather:stability_}
 for precise definitions and proofs).

\hfill

\lemma (\cite{_Mather:stability_}, 
Proposition 11.1)  Let $f: Y\to B$ be a smooth mapping of smooth manifolds and $W$ a closed subset of $Y$ admitting 
Whitney stratification, such that $f:W\to B$ is proper. If the restriction of $f$ to each stratum of $W$ 
is a submersion then $W$ is locally trivial over $B$.

\hfill

The idea is that $W$ acquires a structure of a stratified set so that $f$ is a ``controlled submersion'',
meaning that one can mimick Ehresmann's construction 
of diffeomorphism between the fibers of a smooth,
proper submersion in this setting.

In our situation, we can clearly stratify $U$ by complex analytic subsets and choose a stratification of 
${\cal Z}_U$ in such a way that the condition of the lemma is satisfied above the strata. We obtain that
the family ${\cal Z}_U$ is locally trivial over a complement to a (lower-dimensional) analytic subset in 
$\Teich_z^{\min}$, over a complement to an analytic subset in that analytic subset, etc.

\subsection{Proof of the main result} 
\label{_Proofs_Subsection_} 

We know that the family ${\cal Z}$ is locally 
trivial on the complement to a union (possibly countable,
but finite in a neighbourhood of any point
in the base) of proper analytic subsets $P\subset \Teich_z^{\min}$.
First we pick a point $x\in \Teich_z^{\min}$ which is not
in $P$ and whose $\Gamma$-orbit is dense.
Then $x$ has a neighbourhood $U_x$ over which all fibers ${\cal Z}_b$ are homeomorphic. Moreover the union 
$\bigcup_{\gamma\in \Gamma}\gamma(U_x)$ is a dense open subset of $\Teich_z^{\min}$ and all fibers ${\cal Z}_b$
over this union are homeomorphic.

Take another point $x'\in \Teich_z^{\min}$ with dense $\Gamma$-orbit (we call such a complex structure 
{\bf ergodic}, see \cite{_Verbitsky:ergodic_}). Then the orbit of $x'$ hits 
$\bigcup_{\gamma\in \Gamma}\gamma(U_x)$ and therefore ${\cal Z}_{x'}$ is homeomorphic to ${\cal Z}_b$ for 
$b\in P$.

Now take $y\in \Teich_z^{\min}$ such that the corresponding complex structure is not ergodic but does not 
have maximal Picard number either
(``the intermediate orbit'' of \ref{orbits-per}
and \ref{orbits-perz}). 
If ${\cal Z}_{y}$ is not homeomorphic to ${\cal Z}_b$ for 
$b\in P$, the orbit of $y$ should not hit $P$. As $P$ is open, the orbit closure does not hit $P$ either.
As the orbit closure is irreducible it must be
contained in an irreducible complement of the complement
to $P$, but this is an analytic subvariety.
However, the closure of an intermediate orbit is not
contained in a proper analytic set, even locally  
(\ref{_interme_orbit_not_anal_Proposition_}).
\endproof




\subsection{Closing remarks}
\label{_Variants_subsection_}


Basically the same argument proves \ref{_main_incidentn_Theorem_}. Indeed, it suffices to consider, instead of the family 
${\cal Z}$ over $\Teich_z^{\min}$, the families ${\cal B}$ whose fiber over $I\in \Teich_z^{\min}$ is the Barlet space of rational 
curves of cohomology class $z$, or more generally $\lambda z$ (with a fixed rational $\lambda$), and the incidence family ${\cal J}\subset {\cal X}\times_{\Teich_z^{\min}}{\cal B}$. Notice that this is 
not quite the situation of the first Whitney's lemma as the families are not naturally embedded in other families smooth over 
$\Teich_z^{\min}$, but its generalizations to the singular case do exist in the literature (see e.g. \cite{Ver}, Th\'eor\`eme 4.14, 
Corollaire 5.1). One thus gets the topological triviality of both families.

\hfill

To prove \ref{_main_MRC_Theorem_}, notice that by \ref{_main_incidentn_Theorem_}, $h$ respects the MRC fibrations of the
components of $Z_I$ resp. $Z_{I'}$, so that the open parts of the fibers are homeomorphic. Moreover, since $h$ is a stratified 
diffeomorphism, locally at a general point $x$ of a general fiber $F_I$, the map $h$ restricts to a diffeomorphism $h_F: F_I\to F_{I'}$
with some MRC fiber of a component of $Z_{I'}$. But in every holomorphic tangent direction at $x$ there is a rational curve, and it is
sent to a rational curve through its image by $h_F$. Therefore $h_F$ maps the holomorphic tangent space at $x$ into the holomorphic tangent space at $h_F(x)$, and vice versa, in other words, it is holomorphic at $x$.

\hfill

In  \ref{_main_incidentn_Theorem_}
and \ref{_main_homeo_Theorem_}
we prove that the fibers of natural
families associated with rational curves 
are homeomorphic and stratified diffeomorphic.
However, there is a version of Thom-Mather theory
which gives bi-Lipschitz equivalence of the fibers
over open strata of Thom-Mather stratification
(\cite{_Parusinski:Lip_}, \cite{_Parusinski:Lip_stra_}).
Then the same arguments as above
prove that the homeomorphisms constructed
in \ref{_main_incidentn_Theorem_}
and \ref{_main_homeo_Theorem_} are bi-Lipschitz.

\hfill

{\bf Acknowledgements:}
We are grateful to Fedor Bogomolov for pointing out a potential error
in an earlier version of this work, and to Jean-Pierre
Demailly, Patrick Popescu, Lev Birbrair and Daniel Barlet for useful
discussions. We are especially grateful to
Fabrizio Catanese who explained to us the basics
of Thom-Mather theory and gave the relevant 
reference.

\hfill

{

{\small
\noindent {\sc Ekaterina Amerik\\
{\sc Laboratory of Algebraic Geometry,\\
National Research University HSE,\\
Department of Mathematics, 7 Vavilova Str. Moscow, Russia,}\\
\tt  Ekaterina.Amerik@gmail.com}, also: \\
{\sc Universit\'e Paris-11,\\
Laboratoire de Math\'ematiques,\\
Campus d'Orsay, B\^atiment 425, 91405 Orsay, France}

\hfill

\noindent {\sc Misha Verbitsky\\
{\sc Instituto Nacional de Matem\'atica Pura e
              Aplicada (IMPA) \\ Estrada Dona Castorina, 110\\
Jardim Bot\^anico, CEP 22460-320\\
Rio de Janeiro, RJ - Brasil }\\
also:\\
{\sc Laboratory of Algebraic Geometry,\\
National Research University HSE,\\
Department of Mathematics, 7 Vavilova Str. Moscow, Russia,}\\
\tt  verbit@mccme.ru}.
 }
}

\end{document}